\documentstyle[twoside]{iserlespaper}

\newcommand{\qbin}[2]{\left[\begin{array}{c}#1\\#2\end{array}\right]}

\renewcommand{\E}{{\rm e}}

\begin{document}
\thispagestyle{empty}
\begin{center}
{\Huge From Schr\"odinger spectra

\vspace{18pt}
to orthogonal polynomials,

\vspace{18pt}
via a functional equation}

\vspace{32pt}
{\Large Arieh Iserles\footnote{Department of Applied Mathematics and
Theoretical Physics, University of Cambridge, England.}}

\vspace{36pt}
\parbox[t]{130mm}{{\bf Abstract} The main difference between certain
spectral problems for linear Schr\"odinger operators, e.g. the almost
Mathieu equation, and three-term recurrence relations for orthogonal
polynomials is that in the former the index ranges across $\ZZ$ and in
the latter only across $\Zp$. We present a technique that, by a
mixture of Dirichlet and Taylor expansions, translates the almost
Mathieu equation and its generalizations to three term recurrence
relations. This opens up the possibility of exploiting the full power
of the theory of orthogonal polynomials in the analysis of
Schr\"odinger spectra.

Aforementioned three-term recurrence relations share the property
that their coefficients are almost periodic. We generalize a method of
proof, due originally to Jeff Geronimo and Walter van Assche, to
investigate essential support of the Borel measure of associated
orthogonal polynomials, thereby deriving information on the underlying
absolutely continuous spectra of Schr\"odinger operators.}

\end{center}

\vfill
\eject

\section{The almost Mathieu equation and orthogonal polynomials}
\label{sec1}

The point of departure of our analysis is the {\em almost Mathieu
equation\/} (also known as the {\em Harper equation\/}). We seek
$\lambda\in\RR$ and $\{a_n\}_{n\in\ZZ}\in\ell_1[\ZZ]$ that satisfy
\begin{equation}\label{1.1}
a_{n-1}-2\kappa\cos(\alpha n+\beta)a_n+a_{n+1}=\lambda a_n,\qquad n\in\ZZ.
\end{equation}
$\alpha,\beta$ and $\kappa\neq0$ being given real constants.

The almost Mathieu equation features in a number of applications
\cite{Hofstadter} and has been already extensively studied
\cite{Bellissard,Cycon,Last,Simon}. The purpose of our analysis is not to
reveal new features of the spectrum of \R{1.1} {\em per se\/}, since
the latter is quite comprehensively known. Instead, we intend to
demonstrate that the almost Mathieu equation exhibits an intriguing
connection with orthogonal polynomials, a connection that lends itself
to far reaching generalizations.

Let $\omega=\E^{\I\alpha}$, $b_1=\kappa\E^{\I\beta}$ and
$b_2=\kappa\E^{-\I\beta}=\bar{b}_1$. We rewrite \R{1.1} as
\begin{equation}\label{1.2}
(b_1\omega^{2n}+b_2)a_n=\omega^n(a_{n-1}-\lambda a_n+a_{n+1}),\qquad n\in\ZZ
\end{equation}
and consider the {\em Dirichlet expansion\/}
$$y(z)=\sum_{n=-\infty}^\infty a_n \exp\left\{b_1^{\frac12}\omega^n
z\right\},\qquad z\in\CC.$$
The choice of a specific branch of the square root of $b_1$ is
arbitrary. Note that, since $|\omega|=1$, it is easy to demonstrate that
$\{a_n\}_{n\in\ZZ}\in\ell_1[\ZZ]$ implies convergence of the series
for all $z\in\CC$ \cite{Hardy,Iserles1}.

We multiply \R{1.2} by $\exp \left\{ b_1^{\frac12}\omega^n z\right\}$
and sum for $n\in\ZZ$. Since
$$y'(z)=b_1^{\frac12}\sum_{n=-\infty}^\infty a_n\omega^n
\E^{b_1^{\frac12}\omega^n z},\qquad y''(z)=b_1\sum_{n=-\infty}^\infty a_n
\omega^{2n} \E^{b_1^{\frac12}\omega^n z},$$
we readily deduce that $y$ obeys the functional differential equation
\begin{equation}\label{1.3}
y''(z)+b_2y(z)=b_1^{-\frac12} \left\{ \omega^{-1}
y'(\omega^{-1}z) -\lambda y'(z)+\omega y'(\omega z)\right\}.
\end{equation}
The solution of \R{1.3} is determined uniquely by the values of $y(0)$
and $y'(0)$.

Dirichlet expansions have been employed by Gregor\u{\i} Defel
and Stanislav Molchanov \cite{Derfel1} to investigate the simplified
spectral problem
$$a_{n-1}-2\kappa q^n a_n+a_{n-1}=\lambda a_n,\qquad n\in\ZZ,$$
where $\kappa\in\RR\setminus\{0\}$ and $q>0$, and they have
demonstrated that it obeys another functional differential equation.
Both the Derfel--Molchanov equation and \R{1.3} are a generalization
of the {\em pantograph equation\/}, that has been extensively analysed
in \cite{Kato} and \cite{Iserles1}. However, it is important to
emphasize that \R{1.3} has an important feature that sets it apart from
other functional equations of the pantograph type, namely that, unless
$\omega\in\RR$, its evolution makes sense only for complex $z$ and
proceeds along circles of constant $|z|$, emanating from the origin.

Inasmuch as the equation \R{1.3} can be analysed directly, our next
step entails expanding it in Taylor series. Thus, letting
$$y(z)=\sum_{m=0}^\infty \frac{y_m}{m!} z^m,$$
substitution in \R{1.3} readily yields
\begin{equation}\label{1.4}
y_{m+1}=\sqrt{b_1}^{-\frac12} (2\cos\alpha m-\lambda)y_m-b_2
y_{m-1},\qquad m=1,2,\ldots.
\end{equation}
It is beneficial to treat $y_m$ as a function of the spectral
parameter $\lambda$ and to define
$$\tilde{y}_m(t):=b_2^{-m/2} y_m(-(b_1b_2)^{1/2}t),\qquad
m=0,1,\ldots.$$
Brief manipulation affirms that \R{1.4} is equivalent to
\begin{equation}\label{1.5}
\tilde{y}_{m+1}(t)=\left(t+\frac{2}{\kappa}\cos\alpha m\right)
\tilde{y}_m(t)-\tilde{y}_{m-1}(t),\qquad m=0,1,\ldots.
\end{equation}
To specify the solution of \R{1.5} in a unique fashion we need to
choose $\tilde{y}_0$ and $\tilde{y}_1$, which, of course, corresponds
to equipping \R{1.3} with requisite initial conditions. We note in
passing that $t=-\kappa\lambda$, hence real initial values in
\R{1.3} correspond to $\tilde{y}_0,\tilde{y}_1\in\RR$.

Each solution of \R{1.5} is a linear combination of two linearly
independent solutions. Setting $\sigma=2/\kappa$, we let
\begin{eqnarray}
r_{-1}(t)&\equiv&0,\nonumber\\
r_0(t)&\equiv&1,\nonumber\\
r_{m+1}(t)&=&(t+\sigma\cos\alpha m)r_m(t)-r_{m-1}(t),\qquad m=0,1,\ldots,
\label{1.6}
\end{eqnarray}
and
\begin{eqnarray}
s_{-1}(t)&\equiv&0\nonumber\\
s_0(t)&\equiv&1\nonumber\\
s_{m+1}(t)&=&(t+\sigma\cos\alpha(m+1)) s_m(t)-s_{m-1}(t),\qquad
m=0,1,\ldots. \label{1.7}
\end{eqnarray}
It is trivial to verify that $\{r_m\}_{m\in\Zp}$ and
$\{s_{m+1}\}_{m\in\Zp}$ are linearly independent, hence they span all
solutions of \R{1.5}. Moreover, we observe that each $r_m(t)$ and
$s_m(t)$ is an $m$th degree monic polynomial in $t$.
This is a crucial observation, by virtue of the Favard theorem
\cite{Chihara}: Given any three-term recurrence relation of the form
\begin{eqnarray*}
p_{-1}(t)&\equiv&0,\\
p_0(t)&\equiv&1,\\
p_{m+1}(t)&=&(t+c_m)p_m(t)-d_m p_{m-1}(t),\qquad m=0,1,\ldots,
\end{eqnarray*}
the monic polynomial sequence $\{p_m\}_{m\in\Zp}$ is orthogonal with
respect to some Borel measure $\D\varphi$, i.e.\ 
$$\int_{\RR} p_m(t)p_n(t)\D\varphi(t)=0,\qquad m\neq n.$$

The essential support $\Xi$ of $\D\varphi$ is of central importance to the
matter in hand, since it is an easy consequence of, for example, the
$n$th root asymptotics of orthogonal polynomials \cite{Stahl} that, as
long as $\Xi$ is bounded, the sum $\sum_{m=0}^\infty p_m(t)z^m/m!$
converges for $z\in\CC$, whereas convergence fails for $z\neq0$ for
$t\not\in\Xi$. Thus, travelling all the way back from orthogonal
polynomials to functional equations, and hence to the almost Mathieu
equation, we deduce that, subject to the linear transformation
$\lambda=-t/\kappa$, essential supports of the Borel measures
corresponding to \R{1.6} and \R{1.7} result in the essential spectrum
of \R{1.1}.

We note that, in principle, the Favard theorem falls short of
producing a unique measure and it is entirely possible that there
exist many Borel measures that produce an identical set of monic
orthogonal polynomials. This, however, is ruled out by the determinacy
of the underlying Hamburger moment problem and the latter can be
affirmed for both \R{1.6} and \R{1.7} by the Carleman criterion
\cite{Chihara}.

The remainder of this paper is devoted toward the determination of
$\Xi$. In Section \ref{sec2} we demonstrate, by extending a technique
due to Jeff Geronimo and Walter van Assche, that the essential support
-- and, indeed, the underlying Borel measure -- can be specified
explicitly when $\alpha/\pi$ is rational. We prove in essence that
$\D\varphi$ is a linear combination of Chebyshev measures of the
second kind, supported on a set of disjoint intervals.

Various results on irrational $\alpha/\pi$ are reported in Section
\ref{sec3}. Inasmuch as the general form of $\D\varphi$ is currently a
matter for conjecture, we derive a number of results that, beside
being of interest on their own merit, are fully consistent with known
results about the almost Mathieu operator, in particular with the
theorem that the absolutely continuous spectrum of \R{1.1} is (for
irrational $\alpha/\pi$) a Cantor set. 

The equation \R{1.1} has been extensively studied in the past and,
inasmuch as there are few outstanding conjectures, our knowledge of
the spectrum of the almost Mathieu operator is quite comprehensive.
Although the approach of this paper introduces a new perspective,
there is no claim that Sections \ref{sec2}--\ref{sec3} add to the
current state of knowledge of Schr\"odinger spectra. This state of
affairs is remedied in Section \ref{sec4}, where the framework of our
discussion undergoes a far reaching generalization. Firstly, we
demonstrate that general periodic potentials with a finite number of
Fourier harmonics lend themselves to similar analysis, except that,
instead of orthogonal polynomials, the outcome is a generalized
eigenvalue problem for a certain matrix pencil. Secondly, we prove
that a multivariate extension of the almost Mathieu equation can be
`transformed' by our techniques to a problem in (univariate)
orthogonal polynomials.

A possible application of our analysis, and in particular of Section
\ref{sec4}, is numerical computation of the essential spectrum of
\R{1.1} and of its generalizations. We do not pursue this further in
the present paper.

The observation that \R{1.1} is `almost' a three-term recurrence
relation -- only the index range is wrong -- hence that the almost
Mathieu equation might be connected with orthogonal polynomials, is
not new. Most prominently, it had been made in \cite{Ismail}, where it
motivated a very interesting generalization of Chebyshev polynomials.
The most important innovation in the present paper is in a technique
that reduces the index range from $\ZZ$ to $\Zp$ by passing from
Dirichlet to Taylor expansions and which can be extended to cater for
a substantially more general problem.

\section{Orthogonal polynomials with periodic recurrence coefficients}
\label{sec2}

The focus of our attention in this section is the three-term
recurrence
\begin{eqnarray}
p_{-1}(t)&\equiv&0,\nonumber\\
p_0(t)&\equiv&1,\nonumber\\
p_{m+1}(t)&=&(t-\alpha_m)p_m(t)-p_{m-1}(t),\qquad m=0,1,\ldots, \label{2.1}
\end{eqnarray}
where the sequence $\{\alpha_m\}_{m\in\Zp}$ is $K$-periodic,
\begin{equation}\label{2.2}
\alpha_{m+K}=\alpha_m,\qquad m=0,1,\ldots.
\end{equation}
Note that both \R{1.6} and \R{1.7} assume this form when
$\alpha/\pi$ is rational. Our objective is to determine the Borel
measure that renders $\{p_m\}_{m\in\Zp}$ into an orthogonal polynomial
system (OPS).

In \cite{vanAssche} the authors consider the following problem. Let
$\{Q_m\}_{m\in\Zp}$ be an OPS whose measure has an essential support
$\Xi_0\subseteq[-1,1]$ and let $T$ be a given $N$th degree polynomial.
Setting $\Xi=T^{-1}(\Xi_0)$ (the latter set is, generically, a union
of $\leq N$ disjoint intervals), they derive a new OPS, whose Borel measure
is supported in $\Xi$, explicitly in terms of $\{Q_m\}_{m\in\Zp}$.
This construction is intimately related to the discussion of this
section, except that we need, in a manner of speech, to travel in the
opposite direction. As it turns out, the Borel measure associated with
\R{2.1} inhabits a sets of disjoint intervals and we identify it by
choosing an appropriate polynomial transformation $T$.

Let
$$q_n(t):=p_{(n+1)K-1}(t),\qquad n=0,1,\ldots.$$
Note that $q_{-1}\equiv0$ and, moreover, \R{2.1} and \R{2.2} imply
\begin{equation}\label{2.3}
p_{nK}(t)=(t-\alpha_0)q_{n-1}-p_{nK-2}(t),\qquad n=1,2,\ldots.
\end{equation}
We seek polynomials $\alpha_\ell,\beta_\ell$, $\ell=0,1,\ldots,K-1$,
such that
$$p_{nK+\ell}(t)=a_\ell(t)q_{n-1}(t)-b_\ell(t)p_{nK-2}(t),\qquad
\ell=0,1,\ldots,K-1,\quad n=1,2,\ldots.$$
Because of \R{2.3} and the definition of $q_n$, we have
\begin{equation}\label{2.4}
\begin{array}{rclcrcl}
a_{-1}(t) & \equiv & 1, & \qquad\qquad & b_{-1}(t) & \equiv & 0,\\
a_0(t) & = & t-\alpha_0, & \qquad\qquad & b_0(t) & \equiv & 1.
\end{array}
\end{equation}
We next substitute in the recurrence relation \R{2.1} and, by virtue
of \R{2.2}, obtain
\begin{eqnarray*}
p_{nK+\ell+1}(t)&=&(t-\alpha_{\ell+1})p_{nK+\ell}(t)-p_{nK+\ell-1}(t)\\
&=&(t-\alpha_{\ell+1})\{ a_\ell(t)q_{n-1}(t)-b_\ell(t)p_{nK-2}(t)\}\\
&&\quad\mbox{}-\{a_{\ell-1}(t)q_{n-1}(t)-b_{\ell-1}(t)p_{nK-2}(t)\}.
\end{eqnarray*}
Thus, comparing coefficients, we derive the recurrences
\begin{eqnarray}
a_{\ell+1}(t)&=&(t-\alpha_{\ell+1})a_\ell(t)-a_{\ell-1}(t),
\label{2.5}\\
b_{\ell+1}(t)&=&(t-\alpha_{\ell+1})b_\ell(t)-b_{\ell-1}(t),\qquad
\ell=0,1,\ldots,K-2, \label{2.6}
\end{eqnarray}
which, in tandem with \R{2.4}, determine $\{a_\ell,b_\ell\}_{\ell=0}
^{K-1}$.

Let $\ell=K-1$, then
$$p_{nK-2}(t)=\frac{a_{K-1}(t)q_{n-1}(t)-q_n(t)}{b_{K-1}(t)}$$
and, shifting the index,
$$p_{(n+1)K-2}(t)=\frac{a_{K-1}(t)q_n(t)-q_{n+1}(t)}{b_{K-1}(t)}$$
Substituting both  expressions into
$$p_{(n+1)K-2}(t)=a_{K-2}(t)q_{n-1}(t)-b_{K-2}(t)p_{nK-2}(t)$$
yields the recurrence relation
\begin{equation}\label{2.7}
q_{n+1}(t)=(a_{K-1}(t)-b_{K-2}(t))q_n(t)-\Delta_{n-2}(t)q_{n-1}(t),
\end{equation}
where
$$\Delta_\ell(t)=\det\left[\begin{array}{ll}
a_\ell(t) & a_{\ell+1}(t)\\ b_\ell(t) & b_{\ell+1}(t)
\end{array}\right], \qquad \ell=0,1,\ldots,K-1.$$

We multiply \R{2.6} by $a_\ell(t)$, \R{2.5} by $b_\ell(t)$ and
subtract from each other. This readily affirms by induction that 
$$\Delta_\ell(t)=\Delta_{\ell-1}(t)=\cdots=1$$
and \R{2.7} simplifies into
\begin{equation}\label{2.8}
q_{n+1}(t)=(a_{K-1}(t)-b_{K-2}(t))q_n(t)-q_{n-1}(t),\qquad
n=0,1,\ldots.
\end{equation}
Note that \R{2.8} is consistent with $n=0$, since $q_{-1}\equiv0$. To
further simplify the recurrence, we observe that
$q_0(t)=a_{K-1}(t)-b_{K-2}(t)$, hence, letting
$$\tilde{q}_n(x):=\frac{q_n(t)}{q_0(t)},\qquad n=-1,0,\ldots,$$
where $x=q_0(t)$, we obtain the three-term recurrence
\begin{eqnarray*}
\tilde{q}_{-1}(x)&\equiv&0,\\
\tilde{q}_0(x)&\equiv&1,\\
\tilde{q}_{n+1}(x)&=&x\tilde{q}_n(x)-\tilde{q}_{n-1}(x),\qquad n=0,1,\ldots.
\end{eqnarray*}
Thus, each $\tilde{q}_n$ is an $n$th degree monic polynomial and, by
virtue of the Favard theorem, $\{\tilde{q}_n\}_{n\in\Zp}$ is an OPS.
It can be easily identified as a shifted and scaled {\em Chebyshev
polynomial of the second kind\/},
$$\tilde{q}_n(x)=2^n U_n\left(\Frac12 x\right),\qquad
n=0,1,\ldots.$$
We thus deduce that
\begin{equation}\label{2.9}
q_n(t)=2^nq_0(t)U_n\left(\Frac12 q_0(t)\right),\qquad n=0,1,\ldots.
\end{equation}

Before we identify the underlying Borel measure, let us `fill in' the
remaining values of $p_m$. By definition, $p_{nK-1}=q_{n-1}$,
$p_{n(K+1)-1}=q_n$, hence the recurrence \R{2.1} gives
\begin{eqnarray*}
(x-\alpha_1)p_{nK}(t)-p_{nK+1}(t)&=&q_{n-1}(t),\\
-p_{nK+\ell-1}(t)+(t-\alpha_{\ell+1})p_{nK+\ell}(t)-p_{nK+\ell+1}(t)&=&0,
\qquad \ell=1,2,\ldots,K-3,\\
-p_{(n+1)K-3}(t)+(t-\alpha_{K-1})p_{(n+1)K-2}(t)&=&q_n(t).
\end{eqnarray*}
This is a linear system of equations, which we write as
\begin{equation}\label{2.10}
A_{K-1}{\bf p}_n={\bf q}_n,
\end{equation}
where
$$A_m=\left[\begin{array}{cccccc}
t-\alpha_1 & -1 \\
-1 & t-\alpha_2 & -1\\
 & -1 & t-\alpha_3 & -1\\
 & & \ddots & \ddots & \ddots\\
 & & & -1 & t-\alpha_{m-1} & -1\\
 & & & & -1 & t-\alpha_m
	    \end{array}\right],\qquad m=1,2,\ldots,K-1.,$$
$${\bf p}_n=\left[\begin{array}{c}
p_{nK}(t)\\p_{nK+1}(t)\\\vdots\\p_{(n+1)K-3}(t)\\p_{(n+1)K-2}(t)
\end{array}\right] \qquad\mbox{and}\qquad {\bf q}_n =\left[
\begin{array}{c}q_{n-1}(t)\\0\\\vdots\\0\\q_n(t)\end{array}\right].$$
We expand the determinant of $A_m$ in its bottom row and rightmost
column. This results in a three-term recurrence relation and
comparison with \R{2.4} and \R{2.6} affirms that $\det A_m=b_m(t)$.
Hence, solving \R{2.10} with Cramer's rule, we deduce that there exist
$(K-2)$-degree polynomials $\tilde{a}_\ell$ and $\tilde{b}_\ell$,
$\ell=0,1,\ldots,K-2$, such that
\begin{equation}\label{2.11}
p_{nK+\ell}(t)=\frac{\tilde{a}_\ell(t)
q_{n-1}(t)+\tilde{b}_\ell(t)q_n(t)}{b_{K-1}(t)},\qquad
\ell=0,1,\ldots,K-2.
\end{equation}
Bearing in mind the definition of $a_\ell$ and $b_\ell$,
$$p_{nK+\ell}(t)=a_\ell(t)q_{n-1}(t)-b_\ell(t)p_{nK-2}(t),$$
we obtain from \R{2.11} the identity
$$\frac{\tilde{a}_\ell(t)q_{n-1}(t)+\tilde{b}_\ell(t)q_n(t)}
{b_{K-1}(t)}=a_\ell(t)q_{n-1}(t)-\frac{b_\ell(t)(\tilde{a}_{K-2}(t)
q_{n-2}(t)+\tilde{b}_{K-2}(t)q_{n-1}(t))}{b_{K-1}(t)}.$$
We next substitute
$$q_{n-2}(t)=(a_{K-1}(t)-b_{K-2}(t))q_{n-1}(t)-q_n(t)$$
({\em pace\/} \R{2.8}) and rearrange terms, whereby
\begin{eqnarray*}
&&\{\tilde{b}_\ell(t)-b_\ell(t)\tilde{a}_{K-2}(t)\}q_n(t)\\
&=&\{-\tilde{a}_\ell+a_\ell(t)b_{K-1}(t)-(a_{K-1}(t)b_\ell(t)-b_\ell(t)
b_{K-2}(t)) \tilde{a}_{K-2}(t)-b_\ell(t)\tilde{b}_{K-2}(t)\} q_{n-1}(t).
\end{eqnarray*}
However, consecutive orthogonal polynomials $q_{n-1}$ and $q_n$ cannot
share zeros \cite{Chihara}, therefore both sides of the last equality
identically vanish and we derive the explicit expressions
\begin{eqnarray}
\tilde{a}_\ell(t)&=&a_\ell(t)b_{K-1}(t)-(a_{K-1}(t)b_\ell(t)-b_\ell(t)
b_{K-2}(t))
\tilde{a}_{K-2}(t)-b_\ell(t)\tilde{b}_{K-2}(t),\qquad\label{2.12}\\
\tilde{b}_\ell(t)&=&b_\ell(t)\tilde{a}_{K-2}(t). \label{2.13}
\end{eqnarray}
Letting $\ell=K-2$ in \R{2.13} gives
$\tilde{b}_{K-2}(t)=b_{K-2}(t)\tilde{a}_{K-2}(t)$ and we substitute
this into \R{2.12}. The outcome is
\begin{equation}\label{2.14}
\tilde{a}_\ell(t)=a_\ell(t)b_{K-1}(t)-a_{K-1}(t)b_\ell(t)
\tilde{a}_{K-2}(t).
\end{equation}
In particular, $\ell=K-2$ and the definition of $\Delta_m$ result in
\begin{eqnarray*}
(1+a_{K-1}(t)b_{K-2}(t))\tilde{a}_{K-2}(t)&=&a_{K-2}(t)b_{K-1}(t)
=a_{K-1}(t)b_{K-2}(t)+\Delta_{K-2}(t)\\
&=&1+a_{K-1}(t)b_{K-2}(t).
\end{eqnarray*}
Since $a_{K-1}b_{K-2}\not\equiv-1$, we conclude that
$\tilde{a}_{K-2}\equiv1$ and substitution in \R{2.13} and \R{2.14}
yields the explicit formulae
$$\tilde{a}_\ell=\det\left[\begin{array}{cc}
a_\ell(t) & a_{K-1}(t) \\ b_\ell(t) & b_{K-1}(t)
\end{array}\right],\qquad \tilde{b}_\ell(t)=b_\ell(t),\qquad
\ell=0,1,\ldots,K-2.$$

\vspace{8pt}
\noindent {\bf Theorem 1} The OPS $\{p_m\}_{m\in\Zp}$ has an explicit
representation in the form
\begin{equation}\label{2.15}
p_{nK+\ell}(t)=\frac{1}{b_{K-1}(t)} \left\{
\det\left[\begin{array}{cc} a_\ell(t) & a_{K-1}(t) \\ b_\ell(t) &
b_{K-1}(t) \end{array}\right] q_{n-1}(t)+b_\ell(t)q_n(t)\right\},
\end{equation}
where $n=0,1,\ldots$, $\ell=0,1,\ldots,K-1$ and the OPS
$\{q_n\}_{n\in\Zp}$ satisfies the three-term recurrence \R{2.8}. \QED

\vspace{8pt}
Note that letting $\ell=K-1$ or $\ell=-1$ in \R{2.15}, in tandem with
\R{2.4}, results in $p_{(n+1)K-1}=q_n$ and $p_{nK-1}=q_{n-1}$
respectively, as required.

Let $t\in\Xi$. Then, by the discussion preceding the representation
\R{2.9}, we know that $\frac12 q_0(t)\in[-1,1]$, and there exists
$\theta\in[-\pi,\pi]$ such that $\frac12 q_0(t)=\cos\theta$. Since, by
the definition of Chebyshev polynomials of the second kind,
$$U_n(\cos\theta)=\frac{\sin(n+1)\theta}{\sin\theta},$$
\R{2.9} implies that
$$q_n(t)=2^n q_0(t)\frac{\sin(n+1)\theta}{\sin\theta}.$$
Substitution into \R{2.15} results in
$$p_{nK+\ell}(t)=\frac{2^nq_0(t)}{b_{K-1}(t)\sin\theta} \left\{
\det\left[\begin{array}{cc} a_\ell(t) & a_{K-1}(t) \\ b_\ell(t) &
b_{K-1}(t) \end{array}\right] \sin n\theta+b_\ell(t)\sin(n+1)\theta
\right\}.$$

We next proceed to determine the essential support $\Xi$ of the Borel
measure $\D\varphi$ which corresponds to the OPS $\{p_n\}_{n\in\Zp}$.
As we have already mentioned, this is very similar to the construction
of Geronimo and van Assche in \cite{vanAssche}, with
$\{q_n\}_{n\in\Zp}$, $\frac12 q_0$ and $[-1,1]$ plying the role of
$\{Q_n\}_{n\in\Zp}$, $T$ and $\Xi_0$ respectively. This similarity
notwithstanding, there are some important differences -- not least
that our argument advances in an opposite direction to that of
Geronimo and van Assche -- and we present here a complete derivation
of $\Xi$. We commence by observing that, by virtue of Theorem~1,
everything depends on the support of $\tilde{q}_n(x(t))=2^n
U_n\left(\frac12 q_0(t)\right)$, $n=0,1,\ldots$. Thus, we seek a Borel
measure $\D\psi$ such that
$$I_{n,m}=\int_{-\infty}^\infty \tilde{q}_n(x(t))\tilde{q}_m(x(t))
\D\psi(t)=0,\qquad n,m=0,1,\ldots, \quad n\neq m.$$
It follows from the defintion of $\tilde{q}_n$ that
$$I_{n,m}=2^{n+k}\int_{-\infty}^\infty U_n\left(\Frac12 q_0(t)\right)
U_m\left(\Frac12 q_0(t)\right) \D\psi(t).$$
Similarly to \cite{vanAssche}, we seek the inverse function to
$x=\frac12 q_0(t)$. Let $\xi_1<\xi_2<\cdots<\xi_s$ be all the minima
and maxima of $q_0$ in $\RR$ (of course, $s\leq K-2$) and define
$\xi_0=-\infty$, $\xi_{s+1}=\infty$. In each interval
$[\xi_j,\xi_{j+1}]$, $j=0,1,\ldots,s$, the function $\frac12 q_0(t)$
is monotone, hence it possesses there a well-defined inverse. We
denote it by $X_j(x)$, hence $\frac12 q_0(X_j(x))=x$. Changing the
integration variable, we have
\begin{eqnarray}
I_{n,m}&=&2^{n+m}\sum_{j=0}^s \int_{\xi_j}^{\xi_{j+1}}
U_n\left(\Frac12 q_0(t)\right) U_m\left(\Frac12 q_0(t)\right)
\D\psi(t) \nonumber\\
&=&2^{n+m} \sum_{j=0}^s \int_{\frac12 q_0(\xi_j)}^{\frac12
q_0(\xi_{j+1})} U_n(x)U_m(x)\D\psi(X_j(x)). \label{2.16}
\end{eqnarray}

We recall that $\{U_n\}_{n\in\Zp}$ is an OPS with respect to the Borel
measure $(1-x^2)^{\frac12}\D x$, supported by $x\in[-1,1]$. Thus, for
every $j=0,1,\ldots,s$ we distinguish among the following cases:

\vspace{6pt}
\noindent{\bf Case 1:} $q_0(\xi_j)\leq-2$ and $2\leq
q_0(\xi_{j+1})$.\\
We stipulate that $\D\psi(X_j(x))$ vanishes for all 
$$x\in\left[\Frac12 q_0(\xi_j),\Frac12 q_0(\xi_{j+1})\right]
\setminus[-1,1].$$
Since $X_j$ increases monotonically in $[\xi_j,\xi_{j+1}]$, the
contribution of this interval to \R{2.16} is
\begin{equation}\label{2.17}
2^{n+m}\int_{-1}^1 U_n(x)U_m(x)\D\psi(|X_j(x)|).
\end{equation}

\vspace{6pt}
\noindent{\bf Case 2:} $q_0(\xi_{j+1})\leq-2$ and $2\leq
q_0(\xi_j)$.\\
Likewise, we require that the support of $\D\psi(X_j)$ is restricted
to $[-1,1]$. $q_0$ decreases monotonically within $[\xi_j,\xi_{j+1}]$
and straightforward manipulation affirms that \R{2.17} represents the
contribution of this interval to \R{2.16}.

\vspace{6pt}
\noindent{\bf Case 3:} $\min\{q_0(\xi_j),q_0(\xi_{j+1})\}>-2$ or
$\max\{q_0(\xi_j),q_0(\xi_{j+1})\}<2$.\\
In that case we cannot fit $[-1,1]$ into $[\xi_j,\xi_{j+1}]$, hence we
stipulate that $\D\psi(X_j)$ is not supported in $[\xi_j,\xi_{j+1}]$.

\vspace{6pt}
Let $\nu_1<\nu_2<\cdots<\nu_r$ be all the indices in
$\{0,1,\ldots,s\}$ such that either Case 1 or Case 2 holds. We require
that $r\geq1$. Then \R{2.16} reduces to
$$I_{n,m}=\int_{-1}^1 U_n(x)U_m(x)\sum_{\ell=1}^r
\D\psi(|X_{\nu_\ell}(x)|).$$
Since the Hamburger moment problem for the Chebyshev measure of the
second kind is determinate, it follows that necessarily
$$\sum_{\ell=1}^r \D\psi(|X_{\nu_\ell}(x)|)=(1-x^2)^{\frac12}\D
x,\qquad x\in[-1,1].$$

\vspace{8pt}
\noindent {\bf Theorem 2} The orthogonality measure corresponding to
the OPS $\{p_m\}_{m\in\Zp}$ is supported by
$$\Xi={\cal I}_1\cup{\cal I}_2\cup\cdots\cup{\cal I}_r,$$
where for each $\ell=1,2,\ldots,r$ ${\cal I}_\ell\subseteq
[\xi_{\nu_\ell}, \xi_{\nu_{\ell+1}}]$ is the unique interval such that
$|q_0(t)|=2$ at its endpoints.

\vspace{6pt}
{\bf Proof.} Follows at once from our construction. \QED

\vspace{8pt}
Figure 1 displays two examples of the present construction, for
different cases of $q_0$. In each case $\Xi$ is the union of the
`thick' intervals.

Harking back to \R{1.6} and \R{1.7}, we let
$\alpha_m=-\sigma\cos\alpha m$ and $\alpha_m=-\sigma\cos(m+1)\alpha$,
$m=0,1,\ldots$, respectively, where $\alpha=2\pi L/K$. Thus,
$\{\alpha_m\}_{m\in\Zp}$ is indeed $K$-periodic. Unsurprisingly, the
outcome of our analysis are the familiar spectral bounds
\cite{Bellissard,Cycon}. The merits of our approach are, however, not
just in providing an alternative proof of known results but also in
extending the framework to the multivariate case in Section \ref{sec4}.

We mention in passing that the analysis of this section can be easily
extended to recurrences of the form
\begin{eqnarray*}
p_{-1}(t)&\equiv&0,\\
p_0(t)&\equiv&1,\nonumber\\
p_{m+1}(t)&=&(t-\alpha_m)p_m(t)-\beta_m p_{m-1}(t),\qquad
m=0,1,\ldots,
\end{eqnarray*}
where both $\{\alpha_m\}_{m\in\Zp}$ and $\{\beta_m\}_{m\in\Zp}$ are
$K$-periodic. This, however, is of little relevance to the theme of
this paper.

\vfill\eject
\includegraphics{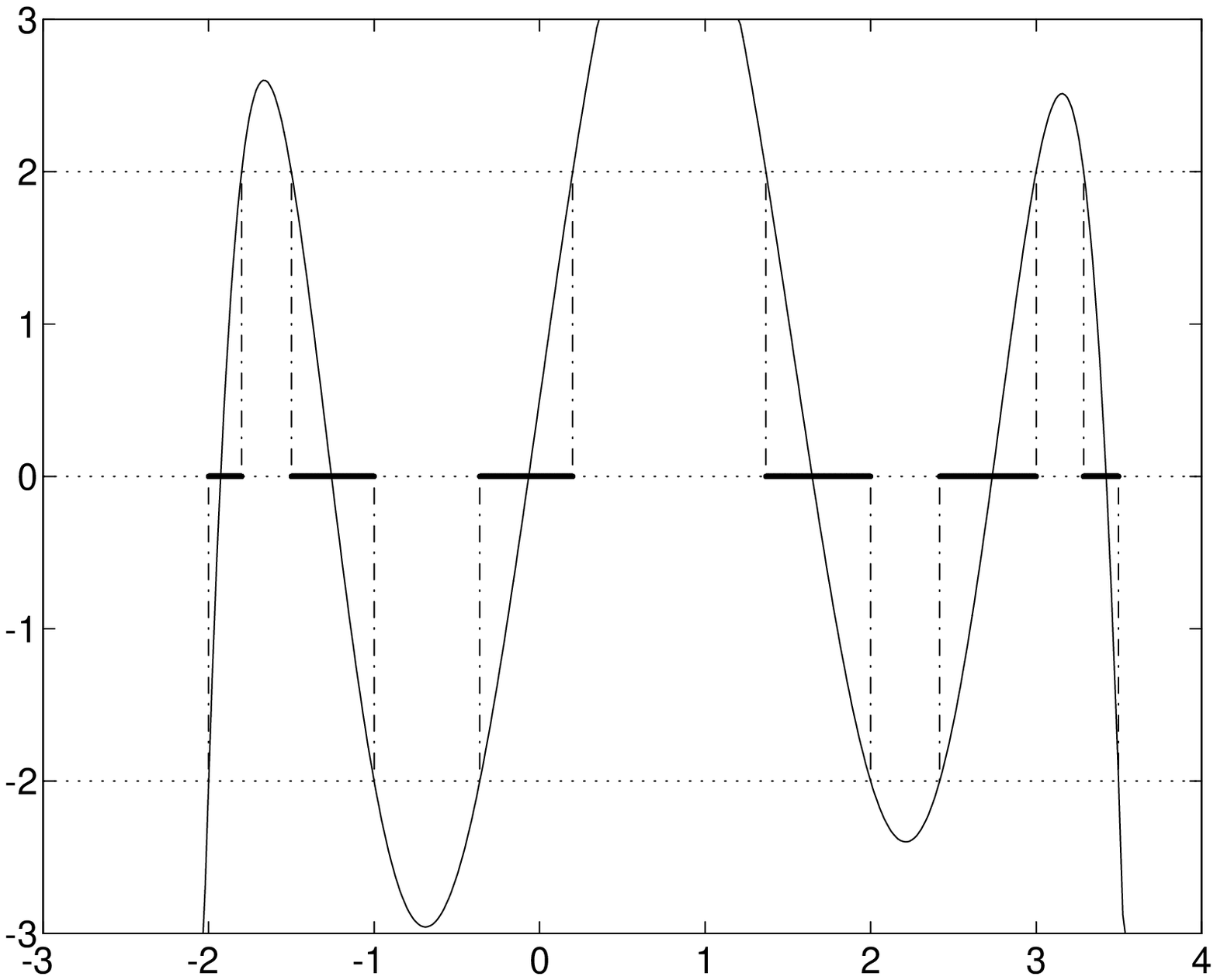}

\includegraphics{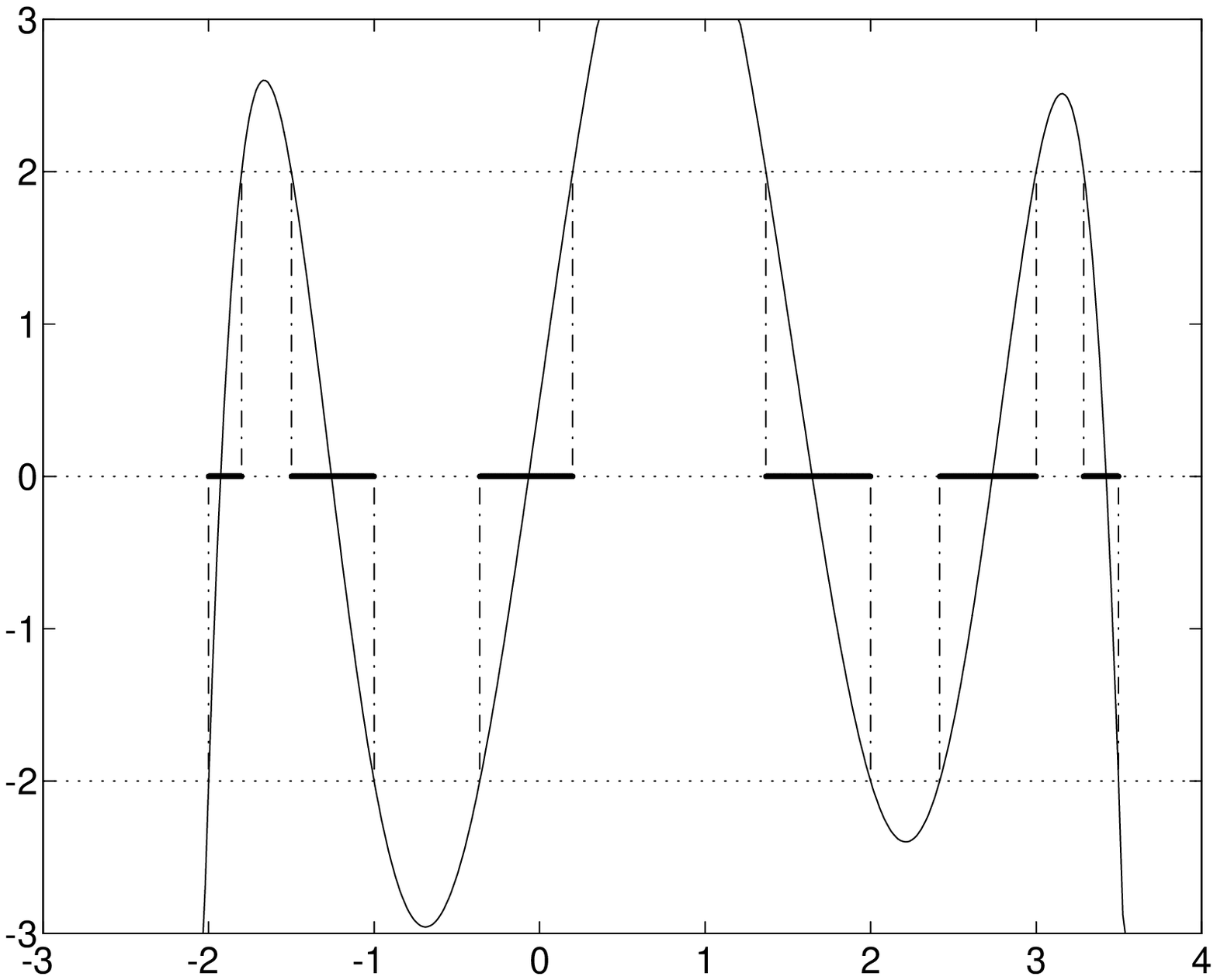}

\rule{0pt}{500pt}

\begin{figure}[b]
\caption{The sets $\Xi$ for two different polynomials $q_0$.}
\end{figure}
\vfill\eject

\section{Orthogonal polynomials with almost periodic recurrence
coefficients} \label{sec3}

The three-term recurrence relations \R{1.6} and \R{1.7} assume,
$\alpha/\pi$ being irrational, almost-periodic recurrence coefficients
and this state of affairs is even more important in a multivariate
generalization of \R{1.1} in Section~\ref{sec4}. Unfortunately, no
general theory exists to cater for orthogonal polynomials with almost
periodic recurrence coefficients. The theme of the present section is
a preliminary and -- in the nature of things -- incomplete
investigation of the case when an irrational $\alpha/\pi$ is
approximated by rationals. In other words, we commence with the
$K$-periodic recurrence \R{2.1}, except that we will allow the period
$K$ to become unbounded.

The motivation for our analysis is an observation which is interesting
on its own merit. Denote by $\sigma_K$ the value of 
$q_{K-1}(0)$ for $\alpha_\ell=-\cos\frac{2\pi\ell}{K}$,
$\ell=0,1,\ldots,K-1$, i.e.\ 
$$\sigma_K=\det\left[\begin{array}{ccccc}
\cos\frac{2\pi}{K} & 1\\
1 & \cos\frac{4\pi}{K} & 1\\
 & \ddots & \ddots\\
 & & 1 & \cos\frac{2(K-2)\pi}{K} & 1\\
 & & & 1 & \cos\frac{2(K-1)\pi}{K} \end{array}\right],\qquad K=1,2,\ldots.$$
Computation indicates that
\begin{equation}\label{3.1}
\sigma_{4L}\equiv0,\qquad
\lim_{L\rightarrow\infty}\sigma_{4L+2}=0,\qquad
\lim_{L\rightarrow\infty} \sigma_{4L+1}=-\lim_{L\rightarrow\infty}
\sigma_{4L+3}=\frac{2\sqrt{3}}{3}.
\end{equation}
It is easy to prove that $\sigma_{4L}=0$ for all $L\geq1$. Thus, let
$K=4L$ and set $\ell=L$ in
\begin{equation}\label{3.2}
\nu_\ell=\cos\frac{2\pi\ell}{K}\nu_{\ell-1}-\nu_{\ell-2}, \qquad
\ell=1,2,\ldots,4L-1
\end{equation}
(note that \R{3.2}, in tandem with $\nu_{-1}=0$, $\nu_0=1$, yields
$\sigma_K=\nu_{K-1}$). This yields $\nu_L+\nu_{L-2}=0$ and we claim
that, in general,
\begin{equation}\label{3.3}
\nu_{L+k}+(-1)^k \nu_{L-k}=0,\qquad k=-1,0,\ldots,L-1.
\end{equation}
We have already proved \R{3.3} for $k=0$ and it is trivially true for
$k=-1$. We continue by induction and assume that \R{3.3} is true for
$k=-1,0,\ldots,s-1$. Letting $\ell=L\pm s$ in \R{3.2}, we have
\begin{eqnarray*}
\nu_{L+s}&=&-\sin\Frac{\pi s}{2L}\nu_{L+s-1}-\nu_{L+s-2},\\
\nu_{L-s}&=&+\sin\Frac{\pi s}{2L}\nu_{L-s-1}-\nu_{L-s-2}.
\end{eqnarray*}
We multiply the second equation by $(-1)^s$ and add to the first, thus
$$(\nu_{L+s}+(-1)^s\nu_{L-s})=-\sin\Frac{\pi s}{2L}
(\nu_{L+s-1}+(-1)^{s-1} \nu_{L-s+1})-(\nu_{L+s-2}+(-1)^{s-2}
\nu_{L-s+2})$$
and \R{3.3} follows at once.

Hence, letting $k=L-1$ in \R{3.3} and recalling that $\nu_{-1}=0$, we
obtain $\nu_{2L-1}=0$. Moreover, similarly to \R{3.3}, we can prove
that 
$$\nu_{3L+k}+(-1)^k\nu_{3L-k}=0,\qquad k=-1,0,\ldots,L-1$$
and $k=L-1$ gives
$$\sigma_{4L}=\nu_{4L-1}=(-1)^L\nu_{2L-1}=0.$$
This completes the proof of the lemma. \QED

\vspace{8pt}
Other observations in \R{3.1} are also true and in the sequel we prove
them in a generalized setting.

Given $\sigma\in(-1,1)$, we define
\begin{eqnarray}
A_{-1}(t)&\equiv&0,\qquad\qquad A_0(t)\equiv1, \nonumber\\
A_n(t)&=&t\xi_n A_{n-1}(t)-A_{n-2}(t),\qquad n=1,2,\ldots, \label{3.4}
\end{eqnarray}
where $\{\xi_n\}_{n=1}^\infty$ is a given real sequence. To emphasize
the dependence on parameters, we write, as and when necessary,
$A_n(\,\cdot\,)=A_n(\,\cdot\,;\xi_1,\xi_2,\ldots,\xi_n)$.

An alternative representation of $A_n$ is 
$$A_n(t)=\det\left[\begin{array}{cccc}
t\xi_1 & 1\\ 1 & t\xi_2 & \ddots\\ & \ddots & \ddots & 1\\ & & 1 &
t\xi_n   \end{array}\right],\qquad n=1,2,\ldots,$$
hence $A_n$ is an $n$th degree polynomial. We observe that
$$A_n(0)=\left\{\begin{array}{lcl}(-1)^s & \qquad &:n=2s,\\0 & \qquad
& :n=2s+1. \end{array}\right.$$
Moreover, differentiating with respect to $t$, we obtain
$$A_n'(t)=\sum_{k=1}^n \det \left[ \begin{array}{ccccccccc}
t\xi_1 & 1 \\
1 & \ddots & \ddots\\
  & \ddots & \ddots & \ddots\\
  & & 1 & t\xi_{k-1} & 1\\
  & & & 0 & \xi_k & 0\\
  & & & & 1 & t\xi_{k+1} & 1\\
  & & & & & \ddots & \ddots & \ddots\\
  & & & & & & \ddots & \ddots & 1\\
  & & & & & & & 1 & t\xi_n   \end{array}\right],$$
therefore, expanding in the $k$th row, we derive the identity
\begin{equation}\label{3.5}
A_n'(t;\xi_1,\ldots,\xi_n) =\sum_{k=1}^n \xi_k
A_{k-1}(t;\xi_1,\ldots,\xi_{k-1}) A_{n-k}(t;\xi_{k+1},\ldots,x_n).
\end{equation}

\vspace{8pt}
\noindent {\bf Proposition 3} $A_n^{(r)}(0)=0$ whenever $n+r$ is odd,
hence the polynomial $A_n$ has the same parity as $n$.

\vspace{6pt}
{\bf Proof.} By induction on $r$. The assertion is true for $r=0$.
Moreover, repeatedly differentiating \R{3.5} with the Leibnitz rule
and letting $t=0$, we obtain
$$A_{2n}^{(2r+1)}(0;\xi_1,\ldots,\xi_{2n})=\sum_{\ell=0}^{2r}
{{2r}\choose \ell} \sum_{k=1}^{2n} \xi_k A_{k-1}^{(\ell)}(0;\xi_1,
\ldots,\xi_{k-1}) A_{2n-k}^{(2r-\ell)}(0;\xi_{k+1},\ldots,\xi_{2n}).$$
But
$$(k-1)+\ell \quad \mbox{is even} \qquad\Longleftrightarrow\qquad
(2n-k)+(2r-\ell) \quad \mbox{is odd},$$
therefore for all $\ell=0,1,\ldots,2r$ and $k=1,2,\ldots,2n$ the
induction hypothesis affirms that at least one of the terms in the
product vanishes. Similar argument demonstrates that
$A_{2n+1}^{(2r)}(0)=0$. \QED 

\vspace{8pt}
We therefore let for all $n=0,1,\ldots$, $s=1,2,\ldots$,
\begin{eqnarray*}
A_{2n}(t;\xi_s,\ldots,\xi_{2n+s-1})&=&\sum_{r=0}^n
B_{2n,s}^{(2r)}t^{2r}\\
A_{2n+1}(t;\xi_s,\ldots,\xi_{2n+s})&=&\sum_{r=0}^n B_{2n+1,s}^{(2r+1)}
t^{2r+1}. 
\end{eqnarray*}
Substitution into \R{3.4} (where we replace $\xi_n$ by $\xi_{n+s-1}$)
results in the recurrences
\begin{eqnarray}
B_{2n,s}^{(2r)}&=&\xi_{s+2n-1}
B_{2n-1,s}^{(2r-1)}-B_{2n-2,s}^{(2r)}, \label{3.6}\\
B_{2n+1,s}^{(2r+1)}&=&\xi_{s+2n}B_{2n,s}^{(2r)}-B_{2n-1,s}^{(2r+1)}.
\label{3.7}
\end{eqnarray}

Given $q\in\CC$, we recall that the {\em $q$-factorial symbol\/} is
defined as
$$(z;q)_n=\prod_{k=0}^{n-1}(1-q^k z),\qquad z\in\CC,\quad
n\in\ZZ\cup\{\infty\},$$
whereas the {\em $q$-binomial\/} reads
$$\qbin{n}{m}:=\frac{(q;q)_n}{(q;q)_m(q;q)_{n-m}},\qquad 0\leq m\leq
n.$$

\vspace{8pt}
\noindent {\bf Lemma 4} Let $\xi_s=q^{\frac12 s}+q^{-\frac12 s}$,
$s=1,2,\ldots$, where $q\in\CC$ is given. Then, for every
$n=0,1,\ldots$, $s=1,2,\ldots$ and $r=0,1,\ldots,n$,
\begin{eqnarray}
B_{2n,s}^{(2r)}&=&(-1)^{n+r}q^{-r\left(2n+s-r-\frac12\right)}
\qbin{n+r}{2r} (-q^{n+s-r};q)_{2r}, \label{3.8}\\
B_{2n+1,s}^{(2r+1)}&=&(-1)^{n+r} q^{-\left(r+\frac12\right) (2n+s-r)}
\qbin{n+r+1}{2r+1} (-q^{n+s-r};q)_{2r+1}.\label{3.9}
\end{eqnarray}

{\bf Proof.} By induction on $n$, using \R{3.8} and \R{3.9}.
Obviously, the assertion of the lemma is true for $n=0$. Otherwise,
for even values,
\begin{eqnarray*}
&&\left(q^{n+\frac{s-1}{2}}+q^{-n-\frac{s-1}{2}}\right)
B_{2n-1}^{(2r-1)}-B_{2n-2,s}^{(2r)}\\
&=&(-1)^{n+r}q^{-\left(r-\frac12\right)(2n+s-r-1)}
\left(q^{n+\frac{s-1}{2}}+q^{-n-\frac{s-1}{2}}\right)
\qbin{n+r-1}{2r-1} (-q^{n+s-r};q)_{2r-1}\\
&&\quad\mbox{}-(-1)^{n+r-1} q^{-r\left(2n+s-r-\frac52\right)}
\qbin{n+r-1}{2r} (-q^{n+s-r-1};q)_{2r}\\
&=&(-1)^{n+r}q^{-r\left(2n+s-r-\frac12\right)} \frac{(q;q)_{n+r-1}}
{(q;q)_{2r}(q;q)_{n-r}} (-q^{n+s-r};q)_{2r-1}\\
&&\quad\mbox{}\times \left\{(1+q^{2n+s-1})(1-q^{2r})
+q^{2r}(1-q^{n-r})(1+q^{n+s-r-1})\right\}\\
&=&(-1)^{n+r}q^{-r\left(2n+s-r-\frac12\right)} \qbin{n+r}{2r}
(-q^{n+s-r};q)_{2r}.
\end{eqnarray*}
This accomplishes a single inductive step for \R{3.8}. We prove
\R{3.9} in an identical manner, by considering odd values of $n$.
\QED

\vspace{8pt}
Recall that our interest in the polynomials $A_n$ has been sparked by
the observation \R{3.1}. Thus, we require to recover cosine terms, and
to this end we choose $q$ of unit modulus.

\vspace{8pt}
\noindent {\bf Proposition 5} Suppose that $q^{n+1}=1$ and $q^m\neq1$ for
$m=1,2,\ldots,n$. Then $B_{2n+1,1}^{(2r+1)}=0$ for all $0\leq r\leq
\frac{n}{2}$.

\vspace{6pt}
{\bf Proof.} Let $2r\leq n$. We have from \R{3.9} that
$$B_{2n+1,1}^{(2r+1)}=(-1)^{n+r}q^{\left(r+\frac12\right)(r+1)}
\frac{(q;q)_{n+r+1}}{(q;q)_{2r+1}(q;q)_{n-r}} (-q^{-r};q)_{2r+1}.$$
However,
$$\frac{(q;q)_{n+r+1}}{(q;q)_{n-r}}=\prod_{\ell=-r}^r (1-q^\ell)=0,$$
whereas, because of our restriction on $r$,
$$(q;q)_{2r+1}=\prod_{\ell=1}^{2r+1}(1-q^\ell)\neq0,$$
since $q$ is a root of unity of {\em minimal\/} degree $n+1$. The
proposition follows. \QED

\vspace{8pt}
\noindent {\bf Corollary} Let $q=\exp\frac{2\pi\I m}{n+1}$, where $m$
and $n$ are relatively prime. Then it is true that
\begin{equation}\label{3.10}
\lim_{n\rightarrow\infty} A_{2n+1}(t)=0
\end{equation}
for every $t\in(-1,1)$.

\vspace{6pt}
{\bf Proof.} Straightforward, since $A_{2n+1}(t)=\O{t^{\frac12 n}}$.
\QED

\vspace{8pt}
\noindent {\bf Proposition 6} Suppose that $\omega=q^{\frac12}$ is a root
of unity of minimal degree $2n+1$. Then, for all $r=0,1,\ldots,n$ it
is true that 
\begin{equation}\label{3.11}
(-1)^n B_{2n,1}^{(2r)}=\prod_{\ell=1}^r \frac{\sin (2\ell-1)\phi}{\sin
2\ell\phi},
\end{equation}
where $\phi=\arg\omega$.

\vspace{6pt}
{\bf Proof.} Since $q^{n+\frac12}=1$, it follows from \R{3.8} that
$$(-1)^n B_{2n,1}^{(2r)}=(-1)^r q^{r\left(r+\frac12\right)}
\frac{(q;q)_{n+r}}{(q;q)_{2r}(q;q)_{n-r}} (-q^{-r+\frac12};q)_{2r}.$$
But
$$\frac{(q;q)_{n+r}(q^{-r+\frac12};q)_{2r}}{(q;q)_{n-r}}=
\prod_{\ell=-r}^{r-1}(1-q^{2\ell+1}).$$
Moreover, $(q;q)_{2r}\neq0$ for $r=0,1,\ldots,n$, since $2n+1$ is the
least nontrivial degree of the root of unity $q$, and we deduce that
$$(-1)^n B_{2n,1}^{(2r)}=(-1)^r q^{r\left(r+\frac12\right)}
\frac{\prod_{\ell=-r}^{r-1}(1-q^{2\ell+1})}{\prod_{\ell=1}^{2r}
(1-q^\ell)}.$$
But
$$\prod_{\ell=1}^{2r}(1-q^\ell)=\prod_{\ell=1}^{2r}\omega^\ell
(\omega^\ell-\omega^{-\ell})=q^{r\left(r+\frac12\right)}
\prod_{\ell=1}^{2r} (\omega^\ell-\omega^{-\ell})$$
and, likewise,
$$\prod_{\ell=-r}^{r-1}(1-q^{2\ell+1})=\prod_{\ell=-r}^{r-1}
\omega^{2\ell+1} (\omega^{2\ell+1}-\omega^{-2\ell-1})=(-1)^r
\prod_{\ell=0}^{r-1} (\omega^{2\ell+1}-\omega^{-2\ell-1})^2.$$
Consequently,
$$(-1)^n B_{2n,1}^{(2r)}=\frac{\prod_{\ell=0}^{r-1} (\omega^{2\ell+1}
-\omega^{-2\ell-1})^2}{\prod_{\ell=1}^{2r}
(\omega^\ell-\omega^{-\ell})}.$$
This is precisely the identity \R{3.11}. \QED

\vspace{8pt}
Next, we consider progression to a limit as $n\rightarrow\infty$ --
\R{3.1} is a special case. Thus, suppose that we have a sequence
$\Phi=\{\phi_n\}_{n\in{\cal I}}$, where $\phi_n=2\pi m_n/(n+1)$,
$m_n\in\Zp$, ${\cal I}\subseteq\Zp$ is a set of infinite cardinality
and
$$\lim_{\stackrel{\scriptstyle n\rightarrow\infty}{n\in{\cal I}}}
\phi_n=\phi\in[0,2\pi).$$
Set
$$C(t,\Phi)=\lim_{\stackrel{\scriptstyle n\rightarrow\infty}{n\in{\cal
I}}} (-1)^{[n/2]}A_n(t;\xi_1^{(n)},\xi_2^{(n)},\ldots,\xi_n^{(n)}),$$
where $\xi_\ell^{(n)}=2\cos\ell\phi_n$, $\ell=1,2,\ldots,n$,
$n=1,2,\ldots$. Thus, $\xi_\ell^{(n)}=q_n^\ell+q_n^{-\ell}$, where
$q_n=\exp\frac{4\pi\I}{n+1}$. Consequently, according to Proposition 5, if
$\cal I$ consists of only odd indices, necessarily $C(t,\Phi)\equiv0$.
This proves, icidentally, that $\sigma_{2L}\rightarrow0$ in \R{3.1}.

\vspace{8pt}
\noindent {\bf Lemma 7} Suppose that ${\cal I}\subseteq 2\Zp$ and that
$\phi=0$. Then, provided that $m_n=o(n^{\frac23})$, it is true that
$C(t,\Phi)=(1-t^2)^{-\frac12}$. 

\vspace{6pt}
{\bf Proof.} Since
$$\frac{\sin(2\ell-1)\phi_n}{\sin2\ell\phi_n}=\frac{2\ell-1}{2\ell}+
\O{\phi_n^3},$$ 
we deduce from \R{3.11} that
$$(-1)^n B_{2n,1}^{(2r)}=r^{-r}{{2r}\choose r}+\O{n\phi_n^3}.$$
According to the assumption, $\O{n\phi_n^3}=o(1)$ and the lemma
follows from
$$\sum_{r=0}^\infty {{2r}\choose
r}\frac{t^{2r}}{4^r}=\frac{1}{\sqrt{1-t^2}}.$$
$\Box$

\vspace{8pt}
Letting $t=\frac12$ affirms the remaining part of \R{3.1}.

Similarly to the last proposition, it is possible to derive an
explicit expression for $C(t,\Phi)$, provided that $\phi/\pi$ is
rational and that $n(\phi-\phi_n)^3=o(1)$ as $n\rightarrow\infty$. The
derivation is long and it will be published elsewhere. It suffices to
mention here the remarkable sensitivity of $C(t,\Phi)$ to both the
choice of $\cal I$ and to the specific nature of $\phi$.

What happens when $\phi/\pi$ is irrational? This is, as things stand,
an open problem. It is possible to show that, formally,
$$C(t,\Phi)=\sum_{\ell=0}^\infty t^{2\ell} \prod_{k=1}^\ell
\frac{\sin(2k-1)\phi}{2k\phi}=
\sum_{\ell=0}^\infty t^{2\ell}q^{\frac14 \ell}
\frac{(q^{\frac12};q)_\ell}{(q;q)_\ell},$$
where $q=\E^{4\I\phi}$. The latter series can be summed up by means of
the Gau\ss--Heine theorem \cite{Gasper} for $|q|<1$ and, after simple
manipulation, for $|q|>1$. Unfortunately, because of a breakdown in
H\"older-continuity across $|q|=1$, it is impossible to deduce its
value on the unit circle from the values within and without by means,
for example, of the Sokhotsky formula \cite{Henrici}.

Clearly, there is much to be done to understand better the behaviour
of the $p_n$'s when the period $K$ becomes infinite. In this section we
have established few results with regard to the values at the origin.
They should be regarded as a preliminary foray into an interesting
problem in orthogonal polynomial theory {\em cum\/} linear algebra and
we hope to return to this theme in the future.

\section{Generalizations} \label{sec4}

There are two natural ways of generalizing an almost Mathieu equation
\R{1.1}, by either specifying a more general periodic potential
or replacing the index by a multi-index. Remarkably, the basic
framework of this paper -- replacing a doubly-infinite recurrence by a
functional equation which, in turn, is replaced by a singly-infinite
recurrence -- survives both generalizations! In the present section we
describe briefly this state of affairs.

Firstly, suppose that the cosine term in \R{1.1} is replaced by a more
general harmonic term and we consider the spectral problem
\begin{equation}\label{4.1}
a_{n-1}-2\left\{\sum_{\ell=1}^m \kappa_\ell \cos(n\ell\theta+\psi_\ell)
\right\} a_n +a_{n+1}=\lambda a_n,\qquad n\in\ZZ.
\end{equation}
We assume that $\kappa_1,\kappa_2,\ldots,\kappa_m\in\RR$ and, without
loss of generality, that $\kappa_m\neq0$. Letting $q=\E^{\I\theta}$,
we set
$$\kappa_\ell^*=\E^{\I\psi_\ell}\kappa_\ell,\quad
\kappa_{-\ell}^*=\E^{-\I\psi_\ell} \kappa_\ell,\qquad
\ell=1,2,\ldots,m,$$
and $\kappa_0^*=0$. Therefore \R{4.1} assumes the form
\begin{equation}\label{4.2}
q^{mn}a_{n-1}-\left\{\sum_{\ell=0}^{2m}
\kappa^*_{\ell-m}q^{n\ell}\right\} a_n+q^{mn}a_{n+1}=\lambda
q^{mn}a_n,\qquad n\in\ZZ.
\end{equation}
Let $c\in\CC\setminus\{0\}$ and consider the Dirichlet series
$y(t)=\sum_{n=-\infty}^\infty a_n \exp\{cq^n t\}$. Since, formally,
$$y^{(\ell)}(t)=c^\ell \sum_{n=-\infty}^\infty q^{n\ell}a_n
\E^{cq^nt},\qquad \ell=0,1,\ldots,$$
we obtain from \R{4.2} the functional differential equation
\begin{equation}\label{4.3}
\sum_{\ell=0}^{2m}\kappa^*_{\ell-m}c^{m-\ell} y^{(\ell)}(t)= q^m
y^{(m)}(qt) -\lambda y^{(m)}(t)+q^{-m}y^{(m)}(q^{-1}t).
\end{equation}
The derivation is identical to that of \R{1.3} and is left to the
reader.

In line with Section \ref{sec1}, we next expand the solution of
\R{4.3} in Taylor series, $y(t)=\sum_{n=0}^\infty p_n t^n/n!$. This
readily yields
$$\sum_{\ell=0}^{2m}\kappa^*_{\ell-m} c^{m-\ell}
p_{n+\ell}=(q^{n+m}+q^{-n-m}-\lambda) p_{n+m},\qquad n=0,1,\ldots,$$
hence, replacing $n+m$ by $n$,
$$\sum_{\ell=-m}^m \kappa_\ell^* c^{-\ell} p_{n+\ell}=(2\cos
n\theta-\lambda) p_n,\qquad \ell=m,m+1,\ldots.$$
Finally, we choose $c=\exp\{\I\psi_m/m\}$, hence
$\kappa_m^*c^{-m}=\kappa^*_{-m}= \kappa_m\in\RR\setminus\{0\}$. We
thus define $\alpha_\ell=c^{-\ell}\kappa_\ell^*/\kappa_m$, $|\ell|\leq
m$ and replace $\lambda$ by $-\lambda \kappa_m$. This results in the
recurrence
\begin{equation}\label{4.4}
\sum_{\ell=-m}^m \alpha_\ell p_{n+\ell}=(\lambda-\beta_n)p_n,\qquad
n=m,m+1,\ldots,
\end{equation}
where
$$\beta_n=-\frac{\cos n\theta}{\kappa_m},\qquad n=m,m+1,\ldots.$$
Note that, inasmuch as the $\alpha_\ell$s may be complex, we have
$\alpha_{-\ell}= \bar{\alpha}_\ell$, $\ell=1,2,\ldots,m$,
$\alpha_0=0$.

The recurrence \R{4.4} is spanned by $2m$ linearly independent
solutions. However, unless $m=1$, it is no longer true that, for
appropriate choice of $p_0,p_1,\ldots,p_{2m-1}$, each $p_n$ is a
polynomial of degree $n+k$ for some $k$, independent of $n$.
Indeed, it is easy to verify that the degree of $p_n$ increases
roughly as $[n/m]$. Hence, orthogonality is lost. Fortunately, an
important feature of orthogonal polynomials, namely that their zeros
are eigenvalues of a truncated Jacobi matrix \cite{Chihara}, can be
generalized to the present framework. It is possible to show that the
zeros of $p_n$ are generalized eigenvalues of a specific pencil of
`truncated' matrices and this provides a handle on their location. We
expect to address ourselves to this issue in a future publication.

Another generalization of \R{1.1} allows the index $n$ to be replaced
by a multi-index ${\bf n}=(n_1,n_2,\ldots,n_d)\in\ZZ^d$. Thus, let
${\bf e}_\ell\in\ZZ^d$ be the $\ell$th unit vector, $\ell=1,2,\ldots,d$,
and consider the spectral problem
\begin{equation}\label{4.5}
\sum_{\ell=1}^d (a_{{\bf n}+{\bf e}_\ell}+a_{{\bf n}-{\bf e}_\ell})
-2\kappa \cos\left(\sum_{\ell=1}^d \alpha_\ell
n_\ell+\beta\right)a_{\bf n}=\lambda a_{\bf n},\qquad {\bf n}\in\ZZ^d.
\end{equation}
In line with Section \ref{sec1}, we let
$$b_1=b\E^{\I\beta},\quad b_2=\E^{-\I\beta},\qquad
q_\ell=\E^{\I\alpha_\ell},\quad \ell=1,2,\ldots,d,$$
whereupon \R{4.5} becomes
\begin{equation}\label{4.6}
\sum_{\ell=1}^d {\bf q}^{\bf n}(a_{{\bf n}+{\bf e}_\ell} +a_{{\bf
n}-{\bf e}_\ell})-(b_1{\bf q}^{2{\bf n}}+b_2)a_{\bf n}=\lambda {\bf
q}^{\bf n}a_{\bf n},\qquad {\bf n}\in\ZZ^d.
\end{equation}
The last formula employs standard multi-index notation, e.g.\ ${\bf
q}^{\bf n}=q_1^{n_1}q_2^{n_2}\cdots q_d^{n_d}$.

We let formally
$$y(t)=\sum_{{\bf n}\in\ZZ^d} a_{\bf n}\exp\left\{ b_1^{\frac12}{\bf
q}^{\bf n}t\right\}$$
and note that
\begin{eqnarray*}
y'(t)&=&b_1^{\frac12} \sum_{{\bf n}\in\ZZ^d} a_{\bf n}{\bf q}^{\bf n}
\exp\left\{ b_1^{\frac12}{\bf q}^{\bf n}t\right\},\\
y''(t)&=&b_1^{\frac12} \sum_{{\bf n}\in\ZZ^d} a_{\bf n}{\bf q}^{2\bf n}
\exp\left\{ b_1^{\frac12}{\bf q}^{\bf n}t\right\}.
\end{eqnarray*}
Therefore, multiplying \R{4.6} by $\exp \left\{ b_1^{\frac12}{\bf
q}^{\bf n}t\right\}$ and summing up for ${\bf n}\in\ZZ^d$ yields,
after brief manipulation, the complex functional differential equation
\begin{equation}\label{4.7}
y''(t)+b_2y(t)=b_1^{-\frac12} \left\{ \sum_{\ell=1}^d
\left(q_\ell^{-1} y'(q_\ell^{-1}t)+q_\ell y'(q_\ell t)\right) -\lambda
y'(t)\right\}. 
\end{equation}

Equation \R{4.7} is of independent interest, being a special case of
the equation
$$y''(t)+c_1y'(t)+c_2 y(t)=\int_0^{2\pi}
y(\E^{\I\theta}t)\D\mu(\theta),$$
where $\D\mu$ is a complex-valued Borel measure. This, in turn, is
similar to the functional integro-differential equations of the form
$$y''(t)+c_1y'(t)+c_2 y(t)=\int_0^1 y(qt)\D\eta(q),$$
say, where $\D\eta$ is, again, a complex-valued Borel measure.
Equations of this kind have been considered by the present
author, jointly with Yunkang Liu \cite{Iserles2}, with an emphasis on
their dynamics and asymptotic behaviour. In the present
paper, however, we are interested in the spectral problem for \R{4.7},
and to this end we again expand $y$ in Taylor series,
$y(t)=\sum_{m=0}^\infty y_m t^m/m!$. It is easy to affirm by
substitution into \R{4.7} the three-term recurrence relation
\begin{equation}\label{4.8}
y_{m+1}=b_1^{-\frac12}\left(2\sum_{\ell=1}^d \cos\alpha_\ell
m-\lambda\right) y_m-b_2y_{m-1},\qquad m=1,2,\ldots.
\end{equation}
Note a most remarkable phenomenon -- although \R{4.5} is
$d$-dimensional, the index in \R{4.8} lives in $\Zp$! In other words,
the dimensionality of the resultant three-term recurrence is
independent of $d$ -- it is, instead, expressed as the number of
harmonics in the recurrence coefficient. Moreover, inasmuch as \R{4.8}
is more complicated for $d\geq2$ then its one-dimensional counterpart
\R{1.4}, both recurrences display similar qualitative characteristics.
In particular, we can use the theory of Section~\ref{sec2} to cater for
the case of $\alpha_1/\pi,\alpha_2/\pi,\ldots,\alpha_d/\pi$ being all
rational. 

In line with the analysis of Section \ref{sec1}, we let
$\tilde{y}_m(t)=b_2^{\frac12 m}y_m(-(b_1b_2)^{\frac12} t)$, $m\in\Zp$,
whereupon \R{4.8} becomes
\begin{equation}\label{4.9}
\tilde{y}_{m+1}(t)=\left(t+\frac{2}{\kappa} \sum_{\ell=1}^d
\cos\alpha_\ell m\right) \tilde{y}_m(t)-\tilde{y}_{m-1}(t),\qquad
m=1,2,\ldots.
\end{equation}
To recover all solutions of \R{4.9} we need to consider a linearly
independent two-dimensional set of solutions. Letting $r_{-1}=0$,
$r_0=1$ and $s_0=0$, $s_1=1$, we recover, similarly to $(1.6$--$7)$,
two sequences $\{r_m\}_{m\in\Zp}$ and $\{s_m\}_{m\in\Zp}$ that span
all solutions of \R{4.9} and such that $\deg r_m=m$, $\deg s_m=m-1$.
In other words, by the Favard theorem both $\{r_m\}_{m\in\Zp}$ and
$\{s_{m+1}\}_{m\in\Zp}$ are OPS and, in line with our analysis of the
one-dimensional almost Mathieu equation \R{1.1}, we are in position
to exploit the theory of orthogonal polynomials.

\section*{Acknowledgements} I have discussed various aspects of this
paper with many colleagues and am delighted to acknowledge their
helpful comments. Particular gratitude deserve Grisha Barenblatt
(Cambridge), Brad Baxter (Manchester), Martin Buhmann (Z\"urich),
Grisha Derfel (Beer Sheva), Mourad Ismail (Tampa), Herb Keller
(Pasadena), Joe Keller (Stanford) and Walter van Assche (Leuven).

\end{document}